\documentclass[12pt]{article}

\usepackage{amssymb,amsmath,amsthm,graphicx}

\newcommand{\eps}{\varepsilon}

\newcommand{\fie}{\varphi}

\newcommand{\abs}[1]{\left|#1\right|}

\newcommand{\bZ}{\mathbb{Z}}

\begin{document}

\begin{center}
{\large Dirichlet's theorem and Jacobsthal's function}
\end{center}

\begin{flushright}
Idris Mercer\\
Florida International University \\
\verb+imercer@fiu.edu+
\end{flushright}

\begin{abstract}
If $a$ and $d$ are relatively prime, we refer to the set of integers
congruent to $a$ mod $d$ as an `eligible' arithmetic progression.
A theorem of Dirichlet says that every eligible arithmetic progression
contains infinitely many primes; the theorem follows from the assertion
that every eligible arithmetic progression contains at least one prime.
The Jacobsthal function $g(n)$ is defined as the smallest positive integer
such that every sequence of $g(n)$ consecutive integers contains an
integer relatively prime to $n$.
In this paper, we show by a combinatorial argument that every eligible
arithmetic progression with $d\le76$ contains at least one prime,
and we show that certain plausible bounds on the Jacobsthal function
of primorials would imply that every eligible arithmetic progression
contains at least one prime.
That is, certain plausible bounds on the Jacobsthal function
would lead to an elementary proof of Dirichlet's theorem.
\end{abstract}

\section{Introduction}

By an {\bf arithmetic progression}, we mean any set of the form
$$
a + d\bZ = \{ a + dn \mid n \in \bZ \}
$$
where $a$ is an integer and $d$ is a positive integer.
For example,
$$
2 + 7\bZ = 9 + 7\bZ = -5 + 7\bZ
= \{ \ldots, -12, -5, 2, 9, 16, \ldots \}.
$$
Note that $a+1\bZ = \bZ$ is an example of an arithmetic progression,
and note that any arithmetic progression $a+d\bZ$ is equal to
$a'+d\bZ$ for some $a'\in\{0,\ldots,d-1\}$.
We will sometimes abbreviate `arithmetic progression' by `AP'.

We say that the arithmetic progression $a+d\bZ$ is {\bf eligible}
if $a$ and $d$ are relatively prime.
This is equivalent to the condition that all elements of $a+d\bZ$
are relatively prime to $d$.

Although our arithmetic progressions contain negative integers,
we will use the word {\bf prime} only for positive primes.
We let $p_k$ denote the $k$th prime, so we have
$$
p_1=2, \qquad p_2=3, \qquad p_3=5, \qquad p_4=7, \qquad p_5=11, \qquad \ldots
$$

Proposition 1.1 below is known to hold for all positive integers $d$,
and the statement that it holds for all $d$ implies Dirichlet's theorem.
The usual proofs of Dirichlet's theorem involve functions of a complex
variable (specifically, Dirichlet L-functions).
Later, we will provide a combinatorial proof that
Proposition 1.1 holds for all $d\le76$, and we will show how certain
plausible bounds on the Jacobsthal function would lead
to a combinatorial proof that Proposition 1.1 holds for all $d$.

{\bf Proposition 1.1.}
If $a+d\bZ$ is an eligible arithmetic progression,
then $a+d\bZ$ contains at least one prime.

{\bf Fact 1.2.}
The statement that Proposition 1.1 holds for {\bf all} positive integers $d$
implies that every eligible arithmetic progression contains {\bf infinitely}
many primes (i.e., Dirichlet's theorem).

{\it Proof.}
Let $a+d\bZ$ be an eligible AP, and by Proposition 1.1,
let $p$ be a prime in $a+d\bZ$. Notice that $a+d\bZ$ is the
disjoint union of $a+2d\bZ$ and $(a+d)+2d\bZ$, which are
eligible APs.
Let $a'+2d\bZ$ be the one of those APs not containing $p$.
Then, applying Proposition 1.1 to $a'+2d\bZ$, we conclude there
exists a prime $p'\in a'+2d\bZ\subset a+d\bZ$,
i.e., there is a prime $p'\ne p$ in $a+d\bZ$.
Next, we can write $a'+2d\bZ$ as the disjoint union
of the two eligible APs $a'+4d\bZ$ and $(a'+2d)+4d\bZ$.
Applying Proposition 1.1 to the one of those APs not containing $p'$,
we can conclude there exists a prime $p''$ distinct from $p$ and $p'$
in a subset of $a+d\bZ$. Continuing in this way, we get an infinite
sequence of distinct primes in $a+d\bZ$. This completes the proof.

We now define primorials and the Jacobsthal function.

{\bf Definition 1.3.}
We define the {\bf primorial} of the $n$th prime to be
$$
p_n\# = \prod_{k=1}^n p_k
$$
so for example, we have $p_1\#=2$, $p_2\#=6$, $p_3\#=30$, $p_4\#=210$, \dots

{\bf Definition 1.4.}
If $n$ is a positive integer, then the
{\bf ordinary Jacobsthal function} $g(n)$ is defined to be the smallest
positive integer $m$ such that among any $m$ consecutive integers,
there is always at least one that is coprime to $n$.

{\bf Example.}
To determine $g(10)$, note that an integer is coprime to 10
if and only if it is congruent to 1, 3, 7, or 9 mod 10.
Thus the longest sequence of consecutive integers that are {\bf not}
coprime to 10 would be a sequence of the form $4+10n,5+10n,6+10n$,
so any sequence of 4 or more consecutive integers must contain at least
one integer that is coprime to 10. Therefore $g(10)=4$.

{\bf Definition 1.5.}
If $n$ is a positive integer, then the
{\bf primorial Jacobsthal function} is defined by
$$
h(n) = g(p_n\#).
$$
So for example, $h(1)=g(2)$, $h(2)=g(2\!\cdot\!3)$,
$h(3) = g(2\!\cdot\!3\!\cdot\!5)$, and so on.

{\bf Example.}
To bound $h(5)=g(2\!\cdot\!3\!\cdot\!5\!\cdot\!7\!\cdot\!11)$,
we first observe that
$$
114, 115, 116, \ldots, 125, 126
$$
are 13 consecutive integers none of which are
coprime to $2\cdot3\cdot5\cdot7\cdot11$.
Next, if we exhaustively check all integers from 1 to
$2\cdot3\cdot5\cdot7\cdot11=2310$,
we find that the above is the longest such sequence,
i.e., any sequence of 14 consecutive integers must contain
an integer coprime to 2310. Thus $h(5)=14$.
(Note that since $p_n\#$ is a rapidly increasing function of $n$,
naive exhaustive search is not necessarily the best way to
calculate $h(n)$ in general.)

Other authors have explored upper and lower bounds for the functions
$g(n)$ and $h(n)$. Such bounds can be expressed either in terms of $n$
or of $p_n$.
An elementary argument shows $h(n) \ge 2p_{n-1}$, and
work of Rankin \cite{Ran}, Maier and Pomerance \cite{MP}, and Pintz \cite{Pin}
leads to lower bounds of the form
$$
h(n) \ge C \cdot \frac{p_n \log p_n \log\log\log p_n}{(\log\log p_n)^2}
$$
(see Section 1 of \cite{Hag}).
As for upper bounds for $h(n)$, Iwaniec \cite{Iwa} showed
$$
h(n) \le C \cdot (n\log n)^2
$$
for an unknown constant $C$. Recalling that the Prime Number Theorem
implies $p_n \sim n\log n$, this is equivalent to $h(n) \le C \cdot p_n^2$.

It has been conjectured that better upper bounds on $h(n)$ are possible.
It was shown in \cite{CW} that
$$
h(n) \le 0.27749612254 \cdot  n^2 \log n
$$
for all $n$ from 50 to 10,000.
Upper bounds smaller than $h(n)=O(p_n^2)$ have interesting consequences.
It was shown in \cite{Kan} that a bound of the form
$$
h(n) \le C \cdot p_n^{2-\eps}
$$
that holds for all $n$ would lead to a short proof of
Linnik's theorem and Dirichlet's theorem.
In Section 2, we will show that if $h(n)$ satisfies the weaker bound
$$
h(n) = o(p_n^2),
$$
for example, if $h(n) = O(p_n^2/\log p_n)$ or $h(n) = O(p_n^2/\log\log p_n)$,
this would lead to a short proof of Dirichlet's theorem.

\section{Main result}

To achieve our main result (Theorem 2.9), it is convenient
to establish some definitions and state some lemmas.

{\bf Definition 2.1.}
A {\bf segment} of an arithmetic progression $X = a+d\bZ$
is any subset of $X$ of the form
$$
\{ a_1, a_1+d, \ldots, a_1+(k-1)d \}
$$
where $a_1 \in X$. We refer to $k$ as the {\bf length}
of the segment.
For example,
$\{ -5, 1, 7, 13 \}$ is a segment of $1+6\bZ$ (of length $4$),
but $\{ -5, 1, 13 \}$ is not a segment of $1+6\bZ$.

{\bf Definition 2.2.}
An {\bf isomorphism} between two arithmetic progressions
$X$ and $Y$ is an order-preserving bijection, i.e.,
a bijection $\fie:X \to Y$ such that $\fie(x)<\fie(x')$
if and only if $x<x'$.

If $\fie$ is a function from $X$ to $Y$, and $Y'\subseteq Y$,
then as is standard, we define
$\fie^{-1}(Y')$ to be the set of all $x \in X$ such that
$\fie(x)\in Y'$.
It is straightforward to verify the following.

{\bf Fact 2.3.}
Suppose $X$ and $Y$ are arithmetic progressions and suppose $\fie$ is an
isomorphism from $X$ to $Y$. If $Y'$ is a segment of $Y$ of length $k$,
then $\fie^{-1}(Y')$ is a segment of $X$ of length $k$.

{\bf Definition 2.4.}
If $c$ is a fixed element of $a+d\bZ$, then we define
the function $\Phi_{c,d}$ from $\bZ$ to $a+d\bZ$ by
$$
\Phi_{c,d}(n) = c + dn.
$$
If the value of $d$ is clear from context, we may write
$\Phi_c = \Phi_{c,d}$.

It is straightforward to verify the following.

{\bf Fact 2.5.}
If $c$ is any fixed element of $a+d\bZ$, then
$\Phi_{c,d}$ is an isomorphism from $\bZ$ to $a+d\bZ$.

{\bf Example.}
Two different isomorphisms from $\bZ$
to $3+5\bZ$ are $\Phi_3 = \Phi_{3,5}$
and $\Phi_{18} = \Phi_{18,5}$.
They are illustrated in the following table.
\begin{center}
\begin{tabular}{ccc}
$n$ & $\Phi_3(n)$ & $\Phi_{18}(n)$ \\ \hline
$-3$ & $-12$ & $3$ \\
$-2$ & $-7$ & $8$ \\
$-1$ & $-2$ & $13$ \\
$0$ & $3$ & $18$ \\
$1$ & $8$ & $23$ \\
$2$ & $13$ & $28$ \\
$3$ & $18$ & $33$
\end{tabular}
\end{center}
Notice that we have, for example,
\begin{align*}
\Phi_3^{-1}(\{3,8,13,18\}) & = \{0,1,2,3\}, \\
\Phi_{18}^{-1}(\{3,8,13,18\}) & = \{-3,-2,-1,0\}.
\end{align*}

{\bf Definition 2.6.}
Let $S$ be a finite set of primes. We say that an isomorphism
$\fie$ from $\bZ$ to $a+d\bZ$ is {\bf $S$-good} if the implication
$$
\mbox{$n$ coprime to all primes in $S$} \Longrightarrow
\mbox{$\fie(n)$ coprime to all primes in $S$}
$$
is true for all $n\in\bZ$.
By a slight abuse of notation, we will sometimes write, e.g.,
$2\!\cdot\!3\!\cdot\!5$-good rather than $\{2,3,5\}$-good.

{\bf Lemma 2.7.}
If $a+d\bZ$ is any eligible arithmetic progression
and $S$ is any finite set of primes, then there is an
$S$-good isomorphism from $\bZ$ to $a+d\bZ$.

{\it Proof.}
Let $q_1,\ldots,q_k$ be all the primes in $S$ not dividing $d$.
Since $d,q_1,\ldots,q_k$ are mutually coprime, the Chinese
Remainder Theorem says there exists an integer $c$ satisfying
all the congruences
\begin{align*}
c & \equiv a \mod d, \\
c & \equiv 0 \mod q_1, \\
& \vdots \\
c & \equiv 0 \mod q_k.
\end{align*}
Note that then $c\in a+d\bZ$ so $c+d\bZ=a+d\bZ$.
We claim that $\Phi_c = \Phi_{c,d}$ is the desired $S$-good isomorphism.
To see this, let $n\in\bZ$, and suppose $n$ is coprime to
all primes in $S$. We must show $\Phi_c(n)=c+dn$ is coprime
to all primes in $S$. Let $p\in S$. If $p\mid d$, then $p\nmid c+dn$
since $c+d\bZ=a+d\bZ$ is an eligible arithmetic progression.
If $p\nmid d$, then $p\mid c$ by construction,
and also $p\nmid n$, so $p\nmid c+dn$. This completes the proof.

{\bf Example.}
If $S=\{2,3\}$, the following are $S$-good
isomorphisms from $\bZ$ to the eligible arithmetic progressions
with $d=3$ or $4$.
\begin{center}
\begin{tabular}{c}
$\Phi_{4,3}$ is a $2\!\cdot\!3$-good isomorphism from $\bZ$ to $1+3\bZ$ \\
$\Phi_{2,3}$ is a $2\!\cdot\!3$-good isomorphism from $\bZ$ to $2+3\bZ$ \\
$\Phi_{9,4}$ is a $2\!\cdot\!3$-good isomorphism from $\bZ$ to $1+4\bZ$ \\
$\Phi_{3,4}$ is a $2\!\cdot\!3$-good isomorphism from $\bZ$ to $3+4\bZ$
\end{tabular}

\begin{tabular}{ccccc}
$n$ & $\Phi_{4,3}(n)$ & $\Phi_{2,3}(n)$ & $\Phi_{9,4}(n)$ & $\Phi_{3,4}(n)$
\\ \hline
$-6$ & $-14$ & $-16$ & $-15$ & $-21$ \\
$\boxed{-5}$ & $\boxed{-11}$ & $\boxed{-13}$ & $\boxed{-11}$ & $\boxed{-17}$ \\
$-4$ & $-8$ & $-10$ & $\boxed{-7}$ & $\boxed{-13}$ \\
$-3$ & $\boxed{-5}$ & $\boxed{-7}$ & $-3$ & $-9$ \\
$-2$ & $-2$ & $-4$ & $\boxed{1}$ & $\boxed{-5}$ \\
$\boxed{-1}$ & $\boxed{1}$ & $\boxed{-1}$ & $\boxed{5}$ & $\boxed{-1}$ \\
$0$ & $4$ & $2$ & $9$ & $3$ \\
$\boxed{1}$ & $\boxed{7}$ & $\boxed{5}$ & $\boxed{13}$ & $\boxed{7}$ \\
$2$ & $10$ & $8$ & $\boxed{17}$ & $\boxed{11}$ \\
$3$ & $\boxed{13}$ & $\boxed{11}$ & $21$ & $15$ \\
$4$ & $16$ & $14$ & $\boxed{25}$ & $\boxed{19}$ \\
$\boxed{5}$ & $\boxed{19}$ & $\boxed{17}$ & $\boxed{29}$ & $\boxed{23}$ \\
$6$ & $22$ & $20$ & $33$ & $27$ \\
$\boxed{7}$ & $\boxed{25}$ & $\boxed{23}$ & $\boxed{37}$ & $\boxed{31}$ \\
$8$ & $28$ & $26$ & $\boxed{41}$ & $\boxed{35}$
\end{tabular}
\end{center}
In this table, there are boxes around all numbers coprime to 6,
illustrating that for each of these isomorphisms $\Phi_{c,d}$,
if $n$ is coprime to all primes in $\{2,3\}$ then also
$\Phi_{c,d}(n)$ is coprime to all primes in $\{2,3\}$,
i.e., each of these isomorphisms is $2\!\cdot\!3$-good.

As another example, consider $S=\{2,3,5\}$ and $a+d\bZ=1+7\bZ$.
Since $120\equiv1$ mod $7$ and $2\!\cdot\!3\!\cdot\!5\mid120$,
the function defined by $\Phi(n)=120+7n$ is a
$2\!\cdot\!3\!\cdot\!5$-good isomorphism from $\bZ$ to $1+7\bZ$.

{\bf Lemma 2.8.}
Suppose $n$ is an integer such that $2 \le n < p_{k+1}^2$
and $n$ is coprime to $p_k\#$. Then $n$ is prime.

{\it Proof.}
Since $n\ge2$, we know $n$ has at least one prime factor.
If $n$ is not prime, then $n$ is a product of two or more primes,
but those primes must be $\ge p_{k+1}$, which would imply $n \ge p_{k+1}^2$.

{\bf Theorem 2.9.}
Let $d$ be a positive integer. If there exists a positive integer $k$
such that
$$
\frac{p_{k+1}^2-2}{h(k)+1} \ge d
$$
then every eligible arithmetic progression $a+d\bZ$
contains at least one prime.

{\it Proof.}
Let $X = \{2,\ldots,p_{k+1}^2-1\}$. Partition $X$ as
$$
X = X_0 \cup \cdots \cup X_{d-1},
$$
where for each $a=0,\ldots,d-1$, we define
$$
X_a = \{ n\in X \mid n \equiv a \pmod d\}.
$$
Then for each $a$, the set $X_a$ is a segment of $a+d\bZ$ of length
$\ge \lfloor(p_{k+1}^2-2)/d\rfloor$.
Now let $a+d\bZ$ be an eligible AP, and let
$S = \{2,3,5,\ldots,p_k\}$.
By Lemma 2.7, there is an $S$-good isomorphism $\Phi$
from $\bZ$ to $a+d\bZ$.
Since $X_a$ is a segment of $a+d\bZ$, we conclude that
$Y=\Phi^{-1}(X_a)$ is a segment of $\bZ$, and we have
$$
\abs{Y} \ge \Big\lfloor\frac{p_{k+1}^2-2}{d}\Big\rfloor \ge
\frac{p_{k+1}^2-2}{d} - 1.
$$
By hypothesis, this implies $\abs{Y} \ge h(k)$.
We conclude that $Y$ contains at least one integer $m$
that is coprime to $p_k\#$. Then, since $\Phi$ is $S$-good,
we conclude that $\Phi(m)\in X_a$ is coprime to $p_k\#$.
Since $\Phi(m) \in X_a \subseteq X = \{2,\ldots,p_{k+1}^2-1\}$,
Lemma 2.8 says $\Phi(m)$ is prime. That is, the eligible AP
$a+d\bZ$ must contain at least one prime.
This completes the proof.

{\bf Corollary 2.10.}
If $d\le76$, then every eligible arithmetic progression $a+d\bZ$ contains
at least one prime. That is, Proposition 1.1 holds for all $d\le76$.

{\it Proof.}
In \cite{ZM}, values of $h(n)$ are computed for $n\le54$, improving upon
\cite{Hag} which computes $h(n)$ for $n\le49$. We find that
$$
\frac{p_{55}^2-2}{h(54)+1} = \frac{257^2-2}{858+1} > 76,
$$
so the result follows from Theorem 2.9.

{\bf Corollary 2.11.}
Suppose $\lim_{n\to\infty} h(n)/p_{n+1}^2 = 0$ (which is the case if,
for example, $h(n) \le C\cdot p_{n+1}^2/\log p_{n+1}$
or $h(n) \le C\cdot p_{n+1}^2/\log\log p_{n+1}$).
Then every eligible arithmetic progression contains at least one prime,
i.e., Proposition 1.1 holds for all $d$.

{\it Proof.}
If $\lim_{n\to\infty} h(n)/p_{n+1}^2 = 0$, then
$\lim_{n\to\infty} p_{n+1}^2/h(n) = \infty$, and then also
$$
\lim_{n\to\infty} \frac{p_{n+1}^2-2}{h(n)+1} = \lim_{n\to\infty}
\frac{p_{n+1}^2}{h(n)}\cdot\frac{p_{n+1}^2-2}{p_{n+1}^2}\cdot\frac{h(n)}{h(n)+1} = \infty
$$
so for every $d$, there exists $k$ with $(p_{k+1}^2-2)/(h(k)+1)\ge d$. The result follows.

{\bf Remark.}
Known values of $h(k)$ appear to support the conjecture that
$(p_{k+1}^2-2)/(h(k)+1)$ grows without bound.
We illustrate with the following table.
\begin{center}
\begin{tabular}{cccc}
$k$ & $p_{k+1}$ & $h(k)$ & $(p_{k+1}^2-2)/(h(k)+1)$ \\ \hline
5 & 13 & 14 & 11.133 \\
10 & 31 & 46 & 20.404 \\
15 & 53 & 100 & 27.792 \\
20 & 73 & 174 & 30.440 \\
25 & 101 & 258 & 39.378 \\
30 & 127 & 330 & 48.722 \\
35 & 151 & 432 & 52.654 \\
40 & 179 & 538 & 59.442 \\
45 & 199 & 642 & 61.585 \\
50 & 233 & 762 & 71.149
\end{tabular}
\end{center}

\bibliographystyle{amsplain}

\end{document}